\documentclass[11pt,a4paper]{article}
\usepackage[hyperref]{acl2021}
\pdfoutput=1
\usepackage{times}
\usepackage{latexsym}

\usepackage{microtype}

\title{Naive probability}
\author{Zal\'an Gyenis\\
  Jagellonian University, Poland\\
  \texttt{zalan.gyenis@gmail.com} \\\And
  Andr\'as Kornai\\
  SZTAKI, Hungary\\
  \texttt{kornai@sztaki.hu} \\}

\date{}

\usepackage{tikz}
\usepackage{latexsym}
\usepackage{amsmath,amssymb,amsfonts,amsthm,amstext,amscd}
\usepackage{mathtools}
\usepackage{multirow}
\usepackage{graphicx}
\usepackage{xcolor}

\usepackage{pgfplots}
\pgfplotsset{compat=1.10}
\usepackage{booktabs}
\usepackage{url}
\usepackage{siunitx}
\usepackage{appendix}
\usepackage{enumerate}
\usepackage{float}

\usepackage{array,makecell}

\aclfinalcopy 

    {\begin{table}[h]\small\centering}%
    {\end{table}}

\begin{document}
\maketitle

\begin{abstract}
We describe a rational, but low resolution model of probability
\end{abstract}

\section{Introduction}

Historically, the theory of probability emerged from the efforts of Pascal and
Fermat in the 1650s to solve problems posed by a gambler, Chevalier de
M\'er\'e \cite{Renyi:1972,Devlin:2008}, and reached its current form in
\cite{Kolmogorov:1933}. Remarkably, not even highly experienced gamblers can
extract high precision probability estimates from observed data: one of de
M\'er\'e's questions concerned comparing the probabilities of getting at least
one 6 in four rolls of one die ($p=0.5177$) and getting at least one double-6
in 24 throws of a pair of dice ($p=0.4914$). Four decades later, Samuel Pepys
is asking Newton to discern the difference between at least two 6s when 12
dice are rolled ($p=0.6187$) and at least 3 6s when 18 dice are rolled
($p=0.5973$).

In this paper we make this phenomenon, the very limited ability of people to
deal with probabilities, the focal point of our inquiry. These limitations, we
will argue, go beyond the well understood limits of numerosity
\cite{Dehaene:1997}, and touch upon areas such as cognitive limits of
deduction \cite{Kracht:2011a} and default inheritance
\cite{Etherington:1987}. We will offer a model of the naive/commonsensical
theory of probability. In Section~\ref{scale} we discuss {\it likeliness},
which we take to be a valuation of propositions on a discrete (seven-point)
scale. In Section~3 we turn to the inference mechanism supported by the naive
theory, akin to Jeffreys-style probability updates. In Section~4 we briefly
sketch the background theory and discuss what we take to be the central
concern, learnability.


\section{The likeliness scale}\label{scale}

We use the term `likeliness' for a valuation on a 7-point scale 0,\ldots,6
which only roughly corresponds to a discretized notion of probability (we
avoid the more natural-sounding `likelihood' as this already has a
well-established technical sense). 0 is assigned to {\it impossible} events,
$l(e)=0$, and 6 to {\it necessary} ones. Note that in this regard $l$
corresponds better to everyday usage in that zero probability events
($p(e)=0$) do occur, and $p(e)=1$ guarantees only that the event $e$ has
measure zero exceptions of occurring. $l(e)=2$ means {\it unlikely}: an
example would be traffic accidents.  $l(e)=1$ means {\it conceivable}, events
that are unlikely in the extreme, but not forbidden by physical law. An
example would be being struck by a meteorite.

There is a duality between $x$ and $6-x$ as in Łukasiewicz $L_7$, so $l(e)=4$
is assigned to {\it likely} events such as travelling without an accident and
$l(e)=5$ to {\it typical} or {\it expected} ones. Almost all lexical knowledge
falls in this last category: chairs are by definition furniture that support a
seated person, and if a particular instance collapses under ordinary weight we
say it failed (whereas we don't conclude that my car failed when I get in a
traffic accident -- alternative hypotheses such as driver error are readily
entertained). Events that are neither likely nor unlikely are assigned the
value 3.

Clearly, using exactly $7$ degrees is somewhat arbitrary, but it is evident
that using only $3$ (say impossible, unknown, possible) would be a gross
oversimplification of how people deal with probability, and using a very fine
scale would create illusory precison that goes beyond people's actual
abilities. With $7$, we stick to a relatively small but descriptive enough
scale. Even if one could argue that, say on cognitive grounds, $5$ or $9$
degrees would be better, the overall methodology would be the exact same, and
everything below could be easily modified and worked out with that scale.
Altogether, our choice of having a $7$-degree scale is more of an illustration
than a commitment, albeit one well supported by practical exerience with
semantic differentials \cite{Osgood:1957}. 

The commonsensical valuation, which is our object of study here, differs from
probabilities in several respects. The most important from our perspective is
{\it lack of additivity}. At this point, it is worth emphasizing that
the theory of likeliness valuation is not intended as a replacement of the
standard (Kolmogorov) notion of probability, which we take to be the correct
theory of the phenomena studied under this heading, but rather as an
explanatory theory of how the {\it naive} worldview accounts for these
phenomena. The fact that as a computational device the standard theory is
superior to the naive theory is no more a reason to abandon study of the naive
theory than the superiority of eukaryotes is reason to abandon study of
prokaryotes.

In this regard, our work differs significantly from studies like
\cite{Pearl:2009} or \cite{Spohn:2012}, which investigate causalty and degrees
of belief in a conservative framework, taking the preservation of all results
of the standard theory for granted. To quote \cite{Spohn:2012} ``Probability
theory is indeed my paradigm; by all means, we must not fall below its
standards.'' Since our goal is to study the commonsensical notion, we are
operating in a different paradigm, and by lack of additivity we don't just
mean lack of $\sigma$-additivity, but something that is already visible on
finite sums. Consider the Law of Total Probability, that $p(A)$ can be
computed as $\sum_n p(A|B_i)p(B_i)$ where the $B_i$ provide a (typically
finite) partition of the event space. The equivalent formulation with
likelinesss normed to 1 would be

\begin{equation}
l(A)=\bigoplus_i l(B_i) \otimes l(B_i \rightarrow A)
\end{equation}

\noindent
Here we retain the assumption that likeliness is a valuation in a semiring
where addition $\oplus$ and multiplication $\otimes$ are defined, but instead
of conditional probability we will speak about relevant implication
$\rightarrow$ having a valuation of its own. The semiring of greatest interest
is the one familiar from $n$-valued logic, where $ \otimes$ is min, and
$\oplus$ is max. We note that we allow two types of propositions only: 
standalone sentences $A$ and sentences in the form of an implication $A\to B$, cf.
Section \ref{sec:3}.

To put lack of additivity in sharp relief, consider the following
commonsensical example: all men are mortal. If we take $A$ to be eventual
death, we have $l(A)=6$. If we ask people to elicit causes of death $B_i$,
they will produce a handful of causes such as cancer or heart attack that they
consider likely ($l=4$); some like accidents of infections they consider
neither very likely nor very unlikely ($l=3$); some like autoimmune diseases
or freezing to death they consider less likely ($l=2$); and some they consider
conceivable but extremely unlikely such as murder/suicide or terrorism
($l=1$). Needless to say, such valuations are not precisely uniform across
people, but they do have high intrasubjective consistency (as measured e.g. by
$\kappa$ statistics). Since $l(B_i \rightarrow A$) is by definition 6, we are
left with an enumeration of causes:

\begin{equation}
l(A)=\bigoplus l(B_i)=\oplus_{i=0}^6 \oplus_{l(B_j)=i} i
\end{equation}

\noindent
The problem here is that no amount of heaping on more of less likely causes
will increase the $\oplus$ above the valuation of its highest term. The
phenomenon is already perceptible at the low end: if we collect all
conceivable causes of death from lightning strike to shark attack, we have
`death by (barely) conceivable causes' which itself is unlikely, not just 
conceivable. 

In actual mortality tables, this phenomenon is reflected in the proliferation
of categories like `unknown', `unspecified', and `other', which take up the
slack. Depending on the depth of tabulation, the catchall category typically
takes up between .5\% and 5\% of the total data, which corresponds well to the 
lack of sensitivity below 1\% observed in the de M\'er\'e and Pepys examples
we started with. 

Another obvious difference between the standard and the naive theory is the
way extremely low or extremely high probability events are treated. When we
want to draw the line between impossible and conceivable events, we don't rely
on a single numerical cutoff. Nevertheless, to get a rough idea of the
probability values corresponding to the various likeliness values, we take the
proverbial `one in a billion chance' as marking, in some fuzzy sense, the
impossible/conceivable boundary.  Using $\log_{10}$ odds scale (for a
justification of base 10, see \cite{Jaynes:2003} Ch. 4.2), this gives $-9$, so
the next natural order of magnitude \cite{Gordon:2017} is $-9/\sqrt{10}
\approx -2.846$. Converting back from log odds to probabilities, this brings
us to $p=0.001423$, which we can take to mark the conceivable/unlikely
boundary. The next natural order of magnitude brings us to log odds $-0.9$,
corresponding to $p=0.1118$, which marks the unlikely/neutral boundary (see
Fig.~\ref{fig1}). 

By symmetry, in this reckoning everything between $p=0.1118$ and $p=0.8882$ is
considered $l=3$, neither particularly likely nor particularly unlikely.
Likely events are between $p=0.8882$ and $p=0.9986$, while typical events are
above that limit though still with a one in a billion chance of failure. As at
the low end, the naive theory lacks the resolution to distinguish such failure
rates from necessity (total absence of failure).

\begin{figure}
\hspace*{-3mm}\resizebox{0.5\textwidth}{!}{%

\begin{tikzpicture}[thick,scale=1, every node/.style={scale=0.9}]
    \begin{axis}
        [
        ,width=17cm
        ,xlabel=p value
        ,ylabel=l value
		,ytick={-1,0, 1, 2, 3, 4, 5, 6}
		,yticklabels={0, 1, 2, 3, 4, 5, 6}
		,y tick label style={yshift=30}
        ,xtick={0.0014, 0.1118, 0.8882, 0.9986 }
        ,xticklabels={\textcolor{red}{0.0014}, \textcolor{red}{0.1118}, \textcolor{red}{0.8882}, \textcolor{red}{0.9986}}
		,extra x ticks={0.0125, 0.2008, 0.7992, 0.9875}
		,extra x tick labels={\textcolor{blue}{0.0125}, \textcolor{blue}{0.2008}, \textcolor{blue}{0.7992}, \textcolor{blue}{0.9875}}
		,every extra x tick/.append style={yshift=-6mm, color = blue}
		,grid=both
        ]
        \addplot+[red, thick,  mark = *,mark options=solid] coordinates
        { (-0.03,-1) (-0.02,0) (0.0014,1) (0.1118,2) (0.8882,3) (0.9986,4) (1.02,5) (1.03,6) };

        \addplot+[blue, dotted, thick, mark size= 1pt, mark = *,mark options=solid] coordinates
        { (-0.03,-1)  (-0.02,0) (0.0125,1) (0.2008,2) (0.7992,3) (0.9875,4) (1.02,5) (1.03,6) };
    \end{axis}
\end{tikzpicture}
}
\caption{$l$ as a function of $p$, starting at $10^{-9}$ {\color{red} (red)}
  and $10^{-6}$ {\color{blue} (blue)}}\label{fig1}
\end{figure}

We should emphasize here that it is the overall logic of the scheme that we
are vested in, not the particular numbers. For example, if we assume an
initial threshold of one in a million instead of one in a billion, the limits
will be at $0.0125$ and $0.2008$ (and by symmetry at $0.7992$ and $0.9875$),
as shown by the blue curve in Fig.~\ref{fig1},
but the major characteristics of the system, such as the `neither likely nor
unlikely' category takes up the bulk of the cases, or that $l=2$ cases are
noticeable, whereas $l=1$ cases are barely detectable, remain unchanged. 

It should be emphasized that such limits, however we set them, are not
intended as a crisp characterization of human classification ability, the
decision boundaries are fuzzy. Returning to lack of additivity, there may well
be several likeley causes of death beyond cancer and heart attack, but no
closed list of such is sufficient for accounting for the fact that eventual
death is necessary. For this, we
need a slack variable that lifts the $\oplus$ of the likely $l=4$ causes to
$l=5$ or $l=6$, which we find in $B_n$ `death by other causes'. We note that
historically old age was seen as a legitimate cause of death, and only very
recently (since the 1980s) do coroner's reports and obituaries find it
necessary to list the failure of a specific organ or subsystem as the cause of
death.

Finally, in contradistinction to the standard theory, $\oplus$ can extend only
to a handful of terms, especially as the terms are implicitely assumed
independent. By the above reckoning, it takes less than 80 unlikely causes to
make one neutral, and less than 8 neutral to make a likely one. The geometry
of the likeliness space is {\it tropical} \cite{Maclagan:2015}, with the
naive theory approximating the log odds (max) semiring. 


\section{Naive inference (likeliness update)}\label{sec:3}

We have two types of propositions: stand alone sentences $A$ and sentences in the
form of an implication $A\to B$. A context is a (finite) collection of propositions, which
can be represented by a directed graph: nodes of the graph denote propositions $A$ and edges
of the graph denote implications $A\to B$. The likeliness function is an evaluation 
acting on the graph: both vertices and edges can have numeric values between $0$ and $6$, 
$0$ representing impossibility, $6$ representing necessity. 

Values $l(A\to B)$ belong to the inner model (for details see Section~4),
therefore they are hardly subject to change.  Take the following example as an
illustration. Snowbird is a ski resort in Utah. Say, for a typical European,
Snowbird is related to travelling, skiing, and snowing with the likelinesss
\def\text{\mbox}

\begin{eqnarray*}
	l( \text{Snowbird}\to\text{travelling}) &=& 5 \\
	l( \text{Snowbird}\to\text{skiing}) &=& 5 \\
	l( \text{Snowbird}\to\text{snowing}) &=& 5 
\end{eqnarray*}

\begin{figure}[b]
\begin{center}
\begin{tikzpicture}
	\node (snowbird) at (0,0) {Snowbird};
	\node (skiing) at (2.3,0) {skiing};
	
	\node (accident) at (1,1.5) {accident};
	\node (skiacc) at (4,1.5) {ski-accident};
	\node (death) at (6.7,1.5) {death};
	
	\node (travelling) at (-1.5,2) {travelling};

	\draw[->] (snowbird) to node[above] {$5$} (skiing);
	\draw[->] (snowbird) to node[left] {$5$} (travelling);
	\draw[->] (snowbird) to node[right] {$1$} (accident);
	\draw[->] (travelling) to node[above] {$4$} (accident);
	\draw[->] (skiing) to node[right] {$0$} (accident);
	\draw[->] (accident) to node[above] {$0$} (skiacc);
	\draw[->] (accident) to[bend left] node[above] {$3$} (death);
	\draw[->] (skiing) to node[above] {$4$} (skiacc);
	\draw[->,dashed] (skiing) to[bend right] node[right] {$3$} (death);
	\draw[->] (skiacc) to node[above] {$3$} (death);
\end{tikzpicture}
\end{center}
\caption{Skiing, accidents, Snowbird}
\end{figure}

Such likelinesses express {\em typicality} of these relations. Skiing is
related to some extent, say, to ski-accident, and ski-accident to death. Take
the example below (for the sake of example we differentiate between
ski-accidents and accidents; the latter excludes accidents occurring while
skiing).

In a typical scenario one does not have any likeliness of the implication
Snowbird$\to$death inside the inner model. However, naive inference works: 
Snowbird {\em typically} implies skiing; skiing is {\em likely} to imply ski-accident; finally,
it is neither likely nor unlikely that ski-accident results in death. Therefore, 
one may say [visiting] Snowbird is neither likely nor unlikely to result in death, i.e.
\begin{eqnarray*}
	l(\text{Snowbird}\to\text{death}) &=& 3
\end{eqnarray*}
In a similar manner, one could obtain the likeliness $l($skiing$\to$death$)=3$ by
saying that skiing is {\em likely} to ensure a ski-accident, while it is neither 
likely nor unlikely that ski-accident results in death.\\

In virtue of the examples above we give a formal model. 
Let assume we have a finite directed graph $G=(V,E)$ and an 
evaluation $l:E\to\{0,\ldots,6\}$. We would like to evaluate edges of the complete
graph on $V$ that are not in $E$. Pick two vertices $a,b\in V$, $a\neq b$ and 
suppose $(a,b)\notin E$. Let $p = (v_1, \ldots, v_n)$ be a path in $G$ from $a=v_1$ 
to $b=v_n$. We write 
\begin{equation}\label{eq:lp}
l(p) = \min\big\{ l(v_i\to v_{i+1}):\; i=1\ldots n-1 \big\} 
\end{equation}

The value $l(p)$ expresses how likely the inference $a\to b$ is in case we are relying
on the chain of already evaluated implications belonging to the path $p$. Then the
value $l(a\to b)$ is obtained as
\begin{equation} 
\min\big\{ l(p):\; p\text{ is a path in $G$ from $a$ to $b$}\big\} 
\end{equation}

In the example above vertices of the graph did not have likelinesss. Suppose
we get new information about John: he is {\em likely} to be in Snowbird, i.e.
$l($Snowbird$)=4$. What consequences can we draw? Being a typical European, if
John is in Snowbird, then he must be travelling and it is really typical that
people travel to Snowbird to ski. The information that $l($Snowbird$)=4$
propagate via the edges of the graph: the likeliness of those propositions
that are related to Snowbird (that is, they are connected by an edge in the
graph to Snowbird) will be updated given new information: $l($travelling$)$
and $l($skiing$)$ become $4$.  In the formal model, given the value $l(a)$ and
a path $p = (v_1,\ldots, v_n)$ from $a=v_1$ to $b=v_n$, using the definition of $l(p)$ 
in equation (\ref{eq:lp}) we can update the likeliness $l(b)$ of $b$ as 
\begin{equation}
  \max\big\{ l(a),l(p):p\text{ is a path from $a$ to $b$} \big\} 
\end{equation}

\noindent
This process of updating iterates: neighbours of just updated vertices get updates in
the next round, etc. Supposing the graph is connected, all vertices are assigned with 
likelinesss as shown of Fig.~\ref{fig2}.

Let us now suppose that we learn that John died abroad. The first column of
Table~1 describes the default likelinesses we assign to various causes of
death, with subsequent columns showing the updates based on whether we learn
($l=6$) that the death took place in Reykjav\'ik, Istanbul, or on a tourist
trip, destination unspecified. Some rows are easy to explain: for example
death at home in bed is considered likely, but if we know that John was on a
tourist trip the implication is that he is not at home, and the likeliness is
demoted to 1. Not 0, because there are extremely unlikely but not
inconceivable scenarios whereby he fell in love with the place, bought a home,
and resettled there, cf. \cite{Jaynes:2003} 5.2.2. This is a scenario that is,
perhaps, worth considering if we know only that John went to Reykyav\'ik or
Istanbul and tourism was merely an inferred, rather than explicitly stated,
goal of the trip, but if we {\it know} it was a tourist trip and nothing more
(last column) this is logically incompatible with being at home.

\begin{figure}
\begin{center}
\begin{tikzpicture}
	\node (snowbird) [label=below:{\textcolor{red}{4}}] at (0,0) {Snowbird};
	\node (skiing) [label=below:{\textcolor{blue}{4}}] at (2.3,0) {skiing};
	
	\node (accident) [label=above:{\textcolor{blue}{4}}] at (1,1.5) {accident};
	\node (skiacc) [label=above:{\textcolor{blue}{4}}] at (4,1.5) {ski-accident};
	\node (death) [label=above:{\textcolor{blue}{3}}] at (6.7,1.5) {death};
	
	\node (travelling) [label=above:{\textcolor{blue}{4}}] at (-1.5,2) {travelling};

	\draw[->] (snowbird) to node[above] {$5$} (skiing);
	\draw[->] (snowbird) to node[left] {$5$} (travelling);
	\draw[->] (snowbird) to node[right] {$1$} (accident);
	\draw[->] (travelling) to node[above] {$4$} (accident);
	\draw[->] (skiing) to node[right] {$0$} (accident);
	\draw[->] (accident) to node[above] {$0$} (skiacc);
	\draw[->] (accident) to[bend left] node[above] {$3$} (death);
	\draw[->] (skiing) to node[above] {$4$} (skiacc);
	\draw[->,dashed] (skiing) to[bend right] node[right] {$3$} (death);
	\draw[->] (skiacc) to node[above] {$3$} (death);
\end{tikzpicture}
\end{center}
\caption{John in Snowbird}\label{fig2}
\end{figure}


\begin{table}\begin{center}
\begin{tabular}{c|c|c|c|c}
Cause of death & Default & Reykjav\'ik & Istanbul & trip  \\ \hline 
in hospital  & 4 &   4     &  5    & 4 \\
by accident (non-ski)& 4 & 4 & 4& 5\\
at home in bed & 4 & 1 &1 &0\\
in war & 1 & 0 &0 & 1\\
by homicide & 1 & 1& 1&1\\
by suicide &2 & 2& 2 & 1\\
by forces of nature & 1 & 4 &1 &2 \\
by ski accident &1 & 2& 1 &1 \\
\end{tabular}
\caption{Likeliness of cause of death}\end{center}
\end{table}

The same logic is operative in the next row (war): since we know there is no
war in Reykjav\'ik or Istanbul the likeliness is demoted to 0, but for a
generic trip it is not, since we do know that there are war zones on the globe
and John may have visited one of these.

We obtain that death by ski accident is less likely in Reykjav\'ik (2) than in
Snowbird (3) not because skiing is inherently more safe in Iceland, but simply
because one can travel to Reykjav\'ik for many reasons, and the likeliness
that one goes skiing there is 3, perhaps 4, whereas to ski in Snowbird is
typical (5). In connection of Reykjav\'ik we are much more likely to think of
death by forces of nature, as there are many natural dangers nearby, from
volcanoes to geysirs and sneaker waves, indeed this class rises to the top
category (4). 

This line also illustrates the nonmonotonic nature of the calculus: in general
we consider death by forces of nature conceivable but unlikely in the extreme
(1), knowing that John went to a tourist trip increases this to 2, but further
learning that he went to Istanbul, not particularly known as a natural danger
zone, demotes this back to 1.

\section{Learning}

As we stated at the outset, our goal is to characterize the {\it naive} view
of probability, an undertaking closely tied to naive physics
\cite{Hayes:1979}, folk psychology \cite{Ravenscroft:2016}, etc. Here we
offer a brief, informal outlook of the entire theory, but deal in detail only
with the question of how the adult system of mental representation, called the
`inner model' above, is formed in regards to probabilities. We introduce
terminology by paraphrase, describing the intended meaning before offering
more formal definitions. Our goal is to stay close to the standard meaning of
these terms, but we do not intend to fully recreate every aspect of the
theories where they originate.

One of the central questions for linguistics, both mathematical/computational
and theoretical, is to characterize how text and meaning are related. The
standard answer, provided by Montague Grammar, is to define a homomorphism
from structured/disambiguated text to formulas of intensional logic. For naive
theory, we use a simpler model composed of (hyper)graphs
\citep{Quillian:1969,Collins:1975}.  For the basic building blocks
of our model, the {\it vertices} of a graph, we assume a large number (about
$10^5$) of ur-objects, roughly one per morpheme or word. In addition to these,
we will have a few technical elements such as the empty node $\cdot$, three
directed connectives `0' (is, isa); `1' (subject); and `2' (object). Our
theory of types is rather skeletal, especially when compared to situation
theory \cite{Barwise:1983,Devlin:1991}, with which we share a great deal of
motivation, especially in regards to common-sense reasoning about real world
situations. When we say that a node is (defeasably) typed as location or
person, this simply means that a 0-edge runs from the node in question to the
{\it location} or {\it person} node.\footnote{An initial implementation is available at
\href{https://github.com/kornai/4lang}{https://github.com/kornai/4lang}, and a
parser translating English and Hungarian text to this style of model is
described in \cite{Recski:2016a,Recski:2018}. For further details on knowledge
representation by hypergraphs see \cite{Kornai:2019}.}

There can be various (n-place) {\bf relations} obtaining between objects but,
importantly, relations can also hold between things {\it construed as}
objects, such as geometrical points with no atomic content, e.g. the corner of
the room {\it is next to} the window', complex motion predicates, e.g. `the
flood {\it caused} the breaking of the dam', and so on. Arguments of relations
will be called {\bf matters}, but they need not be material. We use edges of
type 1 and 2 to indirectly anchor such higher relations, so the subject of
causing will have a 1-edge running from the vertex {\it cause} to the vertex
{\it flood}, and the object, the breaking of the dam, will have a 2-edge
running from {\it cause} to the head of the construction where {\it dam} is
subject of {\it burst}. For ditransitive and higher arity relations, which are
tangential to our main topic here, see \cite{Kornai:2012}.  {\bf Valuations}
are partial mappings from graphs (both from vertices and from edges) to some
small linear order $L$ of scores. There is no analogous `truth assignment'
because in the inner models that are central to the theory, everything is true
by virtue of being present. On occasion we may be able to reason based on
missing signifiers, the dog that didn't bark, but this is atypical and left
for later study.  Here we use $L=\{0,\ldots,6\}$ for probability scores, but
similar scales are standardly used in the measurement and modeling of all
sorts of psychological attitudes since \cite{Osgood:1957}. Of particular
interest is the {\it activity} valuation taking values in $A=\{-1,0,1,2\}$,
where -1 means `blocked', 0 means `inactive', 1 means `active', and 2 means
`spreading'. These are used to keep track of the currently active part of the
graph and implement the {\it spreading activation} model of
\cite{Quillian:1969,Nemeskey:2013}.

Learning, therefore, requires three kinds of processes: the learning of nodes,
the learning of edges, and the learning of valuations. We discuss each in
turn.

\paragraph{Learning new vertices}

We assume a small, inborn set of nodes roughly corresponding to cardinal
points of the body schema \cite{Head:1911} and cardinal aspects of the
outside world such as the gravity vertical \cite{Campos:1970}, to which
further nodes are incrementally adjoined. This typically happens in one shot,
a single exposure to a new object like a {\it boot} is sufficient to set up a
permanent association between the word and the object, likely including
sensory snapshots from smell to texture and a prototypical image
\cite{Rosch:1975}. The association is effected by relations such as spatial
{\it on} `boot on foot' and the more abstract teleological {\it for} `boot for
excursion'. 

On rare occasions, children may learn abstract nodes, such as {\it color},
based on explicit enumerations `red isa color, blue isa color,\ldots', but on
the whole we don't have much use for post hoc taxonomic categories like {\it
  footwear}. Rather, we assume that seeing the boot on a foot, and having
already acquired the notion of {\it shoe}, the child simply adds an edge `boot
isa shoe' to their preexisting representation, a matter to which we now turn. 

\paragraph{Learning new edges}

Again, we assume a small, inborn set of edges (0,1,2), and an inborn mechanism
of spreading activation. The inborn edges are learned by a direct mechanism:
once the edge `boot on foot' is activated, this spreads to nodes associated to
{\it boot} (initially, none) and to {\it foot}. For the sake of the example,
let us assume that the child already knows about shoes. If not, we could start
by describing the earlier learning process, whereby {\it shoe} gets associated
to {\it foot} (which was posited as part of the body schema the child is born
with) without altering our main point, that learning is always incremental
attachment to previously learned nodes. Now, since {\it foot} is activated,
this spreads to {\it shoe}, and the clild adds the new 0-edge between boot and
foot. 

The matter is a bit more complex when the association to be learned is not one
of the primitive ones 0,1,2, but a contentful edge like {\it for}. If the
parents are skinheads, the association `boot for excursion' may never get
formed, since the parents wear the boots on all occasions. But if the boots
are only worn for excursions (or construction work, or any other specific
occasion already identified as such by the child) we will see the {\it boot}
and the {\it excursion} nodes jointly activated, which will prompt the creation
of a new link between the two. 

\paragraph{Learning valuations}

We assume that the activation mechanism is unlearned (innate), and we may
assume that some valuations, in particular the sensory hurt/enjoy valuation,
are at least partially innate, e.g. in regards to harm to body parts. But this
still leaves open the question of how we know that forces of nature are a
likely cause of death in Reykjav\'ik but not in Istanbul? Surely this
knowledge is not innate, and most of us have not studied mortality tables and
statistics at this level of specificity, yet the broad conclusion, that death
by natural forces is more likely in Reykjav\'ik than in Istanbul, is present
in rational thinking at the very least in a defeasible form (we will revise
our naive notions if confronted with strong statistical evidence to the
contrary).

Part of the answer was already provided in Section~3, where we described the
mechanism to compute these values. Aside from very special cases, we assume
that such valuations are always computed afresh, rather than stored. What is
stored are simpler building blocks, such as `volcano near Reykjav\'ik',
`volcano isa danger' from which we can easily obtain `danger near
Reykjav\'ik'. A great deal of background information, such that {\it danger}
is connected to {\it death}, must be pulled in to compute the kind of
valuations we described in Table~1, but this does not alter the main point we
are making here, that inner models are small information objects (the entire
mental lexicon is estimated to be about 1.5MB, see \citet{Mollica:2019}).

\section{Conclusions}

We have offered a rational reconstruction of the naive theory of probability.
This theory is not as powerful computational device as the standard theory,
and generally only leads to rough estimates of likeliness. However, it is
better suited for studying human cognitive behavior, as it requires very
little data, and extends to a broad range of cases where the statistical data
undergirding the standard theory is unavailable.

\section*{Acknowledgments}

Gyenis was supported by the Premium Postdoctoral Grant of the Hungarian
Academy of Sciences and by the Hungarian Scientific Research Found (OTKA),
contract number K115593.

Kornai was supported by the Hungarian Ministry of Innovation and Technology
NRDI Office within the framework of the Artificial Intelligence National
Laboratory Program (MILAB).

\bibliographystyle{acl_natbib}
\bibliography{ratmod}

\end{document}